\documentclass[11pt]{article}
	
	\newcommand{\blind}{0}
	
    \makeatletter
    \renewcommand\section{\@startsection {section}{1}{\z@}%
                                       {-3.5ex \@plus -1ex \@minus -.2ex}%
                                       {2.3ex \@plus.2ex}%
                                       {\normalfont\fontfamily{phv}\fontsize{16}{19}\bfseries}}
    \renewcommand\subsection{\@startsection{subsection}{2}{\z@}%
                                         {-3.25ex\@plus -1ex \@minus -.2ex}%
                                         {1.5ex \@plus .2ex}%
                                         {\normalfont\fontfamily{phv}\fontsize{14}{17}\bfseries}}
    \renewcommand\subsubsection{\@startsection{subsubsection}{3}{\z@}%
                                        {-3.25ex\@plus -1ex \@minus -.2ex}%
                                         {1.5ex \@plus .2ex}%
                                         {\normalfont\normalsize\fontfamily{phv}\fontsize{14}{17}\selectfont}}
    \makeatother
	
	\usepackage{amsmath}
	\usepackage{graphicx}
	\usepackage{enumerate}
	\usepackage{float}
	\usepackage{listings}
	\usepackage{url} 
	
\begin{document}
		
		\def\spacingset#1{\renewcommand{\baselinestretch}%
			{#1}\small\normalsize} \spacingset{1}
		
		\if0\blind
		{
			\title{\bf On the Stratification of Product Portfolios}
			\author{Vikram Govindan, Wei Xie \\ \\
			The Estee Lauder Companies, New York, USA}
			\date{}
			\maketitle
		} \fi
		
		\if1\blind
		{

            \title{\bf \emph{IISE Transactions} \LaTeX \ Template}
			\author{Author information is purposely removed for double-blind review}
			
\bigskip
			\bigskip
			\bigskip
			\begin{center}
				{\LARGE\bf \emph{IISE Transactions} \LaTeX \ Template}
			\end{center}
			\medskip
		} \fi
		\smallskip

	\spacingset{1.5} 

\section*{Abstract}
		
Stratifying commercial product portfolios into multiple classes of decreasing priority, ABCD analysis, is a common supply chain tool. Key planning parameters that drive strategic and execution priorities are tied to the resulting segmentation. These priorities in turn drive supply chain performance. For large product assortments, manual segmentation is infeasible so an automated algorithm is needed. We therefore advocate that careful attention be paid to the design of such an ABCD algorithm and present three key features that can be incorporated into such a calculation to improve its quality \& commercial utility.

\bigskip

	\noindent%
	{\it Keywords:} \emph{ABCD}; SKU classification; SKU stratification; Supply chain portfolio, Herfindahl-Hirschman index, Statistical concentration, Supply chain planning

\section*{Introduction} \label{s:intro}
ABCD analysis is a classic supply chain tool to prioritize product portfolios. Items are grouped into three or four classes of decreasing priority, labeled A through D, typically based on revenue forecasts. The underlying principle of this approach is that the revenue distribution in commercial product portfolios, sorted by item, high to low, tends to follow a power law distribution, with the majority of revenue being driven by a small proportion of items. Items with a higher classification generally receive a higher inventory investment at end nodes in the distribution center network, a higher priority for production, higher run frequencies at manufacturing sites and a greater focus on demand plan quality (mix and completeness), which in turn drives their supply chain performance.

More formally, for a portfolio of $n$ items, with revenue $r_{1} \geq r_{2} \geq ... \geq r_{n}$, the fraction of cumulative revenue to the $k^{th}$ item is
\begin{equation}
C_{k} = \frac{\sum\limits^{k}r_{i}}{\sum\limits^{n}r_{i}} \hspace{20 mm} (where \hspace{2 mm} C_{0} = 0, \hspace{2 mm} C_{k} \leq 1)  \
\end{equation}

For ABCD analysis, we define the function that maps the $k^{th}$ item to its class as

\begin{equation}
    f(k)=
    \begin{cases}
      A & \text{if}\ \hspace{2 mm} (C_{k} \leq t_{a}) \hspace{1 mm} \text{or} \hspace{1 mm} (C_{k-1} < t_{a} < C_{k})  \\
      B & \text{if}\ \hspace{2 mm} (t_{a} < C_{k-1} < C_{k} \leq t_{b}) \hspace{1 mm} \text{or} \hspace{1 mm} (C_{k-1} < t_{b} < C_{k}) \\
      C & \text{if}\ \hspace{2 mm} (t_{b} < C_{k-1} < C_{k} \leq t_{c}) \hspace{1 mm} \text{or} \hspace{1 mm} (C_{k-1} < t_{c} < C_{k}) \\
      D & \text{if}\ \hspace{2 mm}  t_{c} < C_{k-1} < C_{k}
    \end{cases}
  \end{equation}

$t_{a}$, $t_{b}$ and $t_{c}$ represent revenue thresholds for the segmentation. Typical parameter values are $t_{a} = 0.25,\ t_{b} = 0.65,\ t_{c} = 0.95$. The $or$ clause in Eq.[2] handles two special cases:

\begin{enumerate}
  \item Guarantee at least one A item. For example:\\
  $f(1) = A$ if $C_{1} = 0.4$ and $t_{a} = 0.25$
  \item To pull items that straddle thresholds into the higher class. For example:\\
    $f(9) = A$ if $C_{8} = 0.23$ , $C_{9} = 0.26$ and $t_{a} = 0.25$
\end{enumerate}

For a large multi-brand product portfolio, the impact of not handling these edge cases, as illustrated in table 1, can be material in terms of lower inventory investment and potential lost revenue for items (that straddle threshold boundaries) labeled with a lower class.

\begin{table}[htb]
\centering
\begin{tabular}{|l|l|l|l|l|l|} 
\hline
\textbf{$k$}     & \textbf{Item}   & \textbf{Revenue ($r_{i}$)} & \textbf{Cumulative Revenue} & \textbf{$C_{k}$}   & \textbf{$f(k)$}  \\ 
\hline
1              & Item 1          & \$100               & \$100                       & 0.18          & A              \\ 
\hline
\textbf{2}     & \textbf{Item 2} & \textbf{\$90}       & \textbf{\$190}              & \textbf{0.35} & \textbf{A}     \\ 
\hline
3              & Item 3          & \$80                & \$270                       & 0.49          & B              \\ 
\hline
4              & Item 4          & \$70                & \$340                       & 0.62          & B              \\ 
\hline
\textbf{5}     & \textbf{Item 5} & \textbf{\$60}       & \textbf{\$400}              & \textbf{0.73} & \textbf{B}     \\ 
\hline
6              & Item 6          & \$50                & \$450                       & 0.82          & C              \\ 
\hline
7              & Item 7          & \$40                & \$490                       & 0.89          & C              \\ 
\hline
8              & Item 8          & \$30                & \$520                       & 0.95          & C              \\ 
\hline
9              & Item 9          & \$20                & \$540                       & 0.98          & D              \\ 
\hline
10             & Item 10         & \$10                & \$550                       & 1             & D              \\ 
\hline
\textbf{Total} &                 & \textbf{\$550}      &                             &               &                \\
\hline
\end{tabular}
\caption{Item 2 \& Item 5 that straddle cumulative revenue thresholds are pulled into the higher class. $t_{a} = 0.25,\ t_{b} = 0.65,\ t_{c} = 0.95$}
\end{table}

There are several practical considerations when implementing such a calculation in the supply chain, that we will not focus on here. They are however worth mentioning:

\begin{enumerate}
  \item The use of gross profit rather than revenue to drive the algorithm, based on availability of reliable cost data at the required granularity.
  \item The source (typically the consensus demand plan) and horizon (typically the next 12 months) of revenue forecasts ($r_{i}$) used for the calculation.
  \item The frequency at which the segmentation is refreshed (once a month, once a quarter etc.) that trades off stability vs. agility within the supply chain.
  \item The level of granularity of the segmentation. Typically supply chains have at least two ABCD classifications, at the item level globally (global $ABCD$) and at the item/distribution-center (DC) level locally (local $ABCD$). An item could therefore have a single global class, but several local classes based on relative revenue by DC. An item could be a global C, but a local A in Japan for instance, for an item popular in that region. The global class drives production prioritization at manufacturing sites, while the local class drives local inventory policy.
  \item Efficient implementation and automation of the calculation, especially for portfolios with hundreds of thousands of items or item-DC combinations. 
  \item Promotional items, without explicit revenue, still need to be produced and may need an ABCD classification based on units. 
\end{enumerate}

In the following sections, we limit our attention to three aspects of the calculation itself, that we believe lead to a segmentation with greater commercial utility:

\begin{enumerate}
    \item \textbf{Segmentation}: Applying $f(k)$ to different collections of items leveraging the product hierarchy tree, to identify A items. For example, first separating the product portfolio into new products ($<$ 12 months of sales history) and in-line products ($\geq$ 12 months of sales history) and applying $f(k)$ to the items in each of those groups separately, to identify A items. Then for instance, partitioning that same collection of in-line items by category and applying $f(k)$ again to the collection of items within each category to guarantee A item representation for each category. The second pass identifies additional A items within each category, not labeled as A in the first pass.
    \item \textbf{Concentration}: Addressing the lack of commercial utility of a very few items being classified as A's due to revenue concentration, by raising $t_{a}$ based on a quantitative concentration measure such as the Herfindahl-Hirschman index.
    \item \textbf{Productivity}: Some interesting observations on the productivity (revenue per item) function and the ability to find $t_{a}$ that maximizes $blended$ productivity of the $A$ class in the special case where a certain number of items are already assigned to the class.
\end{enumerate}

\section*{Segmentation Along Product Hierarchy Dimensions} \label{s:sec2}
Large supply chain product portfolios, with tens of thousands of items typically have very detailed product hierarchies. An example of such a hierarchy tree is division, brand, major category, category, subcategory and application. Every item falls uniquely within the hierarchy tree (figure 1).

\begin{figure}[htb]
\includegraphics[scale=0.8]{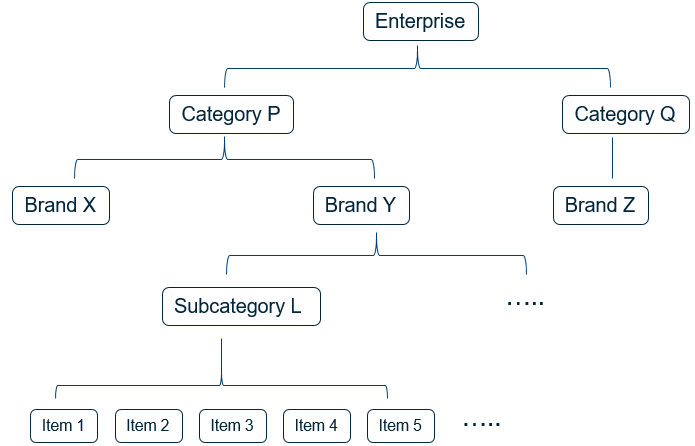}
\caption{An example product hierarchy tree}
\end{figure}

This complexity raises a few challenges:

\begin{enumerate}
    \item Applying $f(k)$ over all items, can result in over representation of some sections of the hierarchy vs. others within the A class. Small but growing brands or categories will be chronically underrepresented in the A class for instance, because of lower revenue items, relative to larger brands or categories. Although technically justifiable, this under representation of A items in smaller brands or categories lacks commercial utility as businesses typically look to drive growth across categories and brands. We formally define under representation as 
    
\begin{equation}
max(C_{k}) < t_{a} \hspace{2 mm} | \hspace{2 mm} \forall \hspace{2 mm} k  \hspace{2 mm} where  \hspace{2 mm}  f(k) = A
\end{equation}
    
    \item New product launches with greater revenue uncertainty and lower productivity relative to  in-line items may merit special consideration due to commercial significance. 
    \item Conversely, for larger brands or categories, items that fall above the cumulative revenue thresholds $t$ (relative to other items within the brand or category) may be significant in absolute dollars and deserve a higher classification.
    \item In the limiting case, applying $f(k)$ for the collection of items at every unique intersection of product hierarchy dimensions could generate too many A's.
\end{enumerate}

To address these challenges, we propose a classification algorithm that proceeds in two stages, the first focused on identifying A items, followed by a second to classify $B$, $C$ and $D$ items. Within the first stage, we propose a series of passes, with each successive pass applying $f(k)$ to identify A items within a specific product hierarchy dimension that was underrepresented over previous passes (figure 2). The choice of dimensions along which to apply $f_{k}$ is driven by commercial relevance, in alignment with key stakeholders.

\begin{figure}[htb]
\includegraphics[scale=0.6]{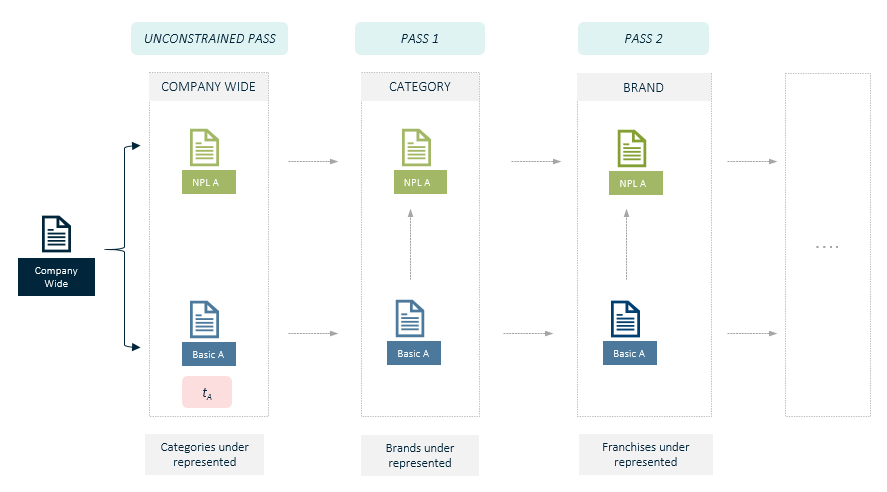}
\caption{Multiple passes to identify A items, first unconstrained, and then along product hierarchy dimensions}
\end{figure}

As an example, in the first stage, the product portfolio is split into new product launches and in-line items and the passes are applied independently. The first pass is unconstrained to capture the largest revenue new and in-line items across the company as A's regardless of position (brand/category/subcategory) in the product hierarchy. Subsequent passes address under representation along different dimensions. For a stable product portfolio, the choice of dimensions and number of passes may be updated once a year. The passes are cumulative along dimensions, so a brand that is over represented in the unconstrained pass, i.e. where $max(C_{k}) >> t_{a} \hspace{2 mm} | \hspace{2 mm} \forall \hspace{2 mm} k  \hspace{2 mm} where  \hspace{2 mm}  f(k) = A$ will receive no additional A's in the subsequent brand pass, but may receive additional A items in passes along other dimensions. See illustrative example in table 2.

\begin{table}[H]
\centering
\begin{tabular}{|l|l|l|l|l|l|} 
\hline
\textbf{Item}   & \textbf{Category} & \textbf{Revenue ($r_{i}$)} & \textbf{Cumulative Revenue} & \textbf{$C_{k}$} & \textbf{$f_{k}$}     \\ 
\hline
Item 1          & P                 & \$1000                     & \$1000                      & 0.05                                                                & A           \\ 
\hline
Item 2          & P                 & \$950                      & \$1950                      & 0.0975                                                              & A           \\ 
\hline
Item 3          & P                 & \$900                      & \$2850                      & 0.1425                                                              & A           \\ 
\hline
Item 4          & P                 & \$850                      & \$3700                      & 0.185                                                               & A           \\ 
\hline
\textbf{Item 5} & \textbf{Q}        & \textbf{\$800}             & \textbf{\$4500}             & \textbf{0.225}                                                      & \textbf{A}  \\ 
\hline
\textbf{Item 6} & \textbf{R}        & \textbf{\$750}             & \textbf{\$5250}             & \textbf{0.2625}                                                     & \textbf{A}  \\
\hline
\end{tabular}
\caption{As an illustration, a first pass across the entire portfolio, identifies items primarily from category $P$ with only one item each from category $Q$ and category $R$. If the revenue for item 5, \$800, represents $< t_{a}$ of the total revenue of category $Q$ as a whole, then category $Q$ is under represented and a second pass applying $f(k)$ on category $Q$ items is warranted for commercial utility. Assumes total portfolio revenue = \$20000}
\end{table}

Once the first stage is complete, and A items have been identified, a second stage to stratify the remaining items into $B$, $C$ and $D$ follows. There is usually not the need for multiple passes in this second stage. After A items have been excluded, if $f(k)$ is applied on the \textit{remainder} of the portfolio, the thresholds for $B$, $C$ and $D$ items are recalculated as

\begin{equation}
t_{b}^{'} = \frac{t_{b} - t_{a}}{1 - t_{a}}
\end{equation}

\section*{Varying Revenue Thresholds Based on Concentration}\label{s:sec3}
So far we have assumed that the cumulative revenue thresholds $t_{a}..t_{c}$ are fixed. We now provide a quantitative framework to adjust these thresholds where revenue is concentrated in a very small proportion of items within a product hierarchy dimension. This can result in very few, sometimes a single A class item, which is not always aligned with commercial go-to-market strategies that endeavor to reliably supply consumers with a wide product assortment.

To get around this difficulty, we propose using the Herfindahl-Hirschman index, a statistical measure of concentration to characterize the product portfolio within relevant product hierarchy dimensions. The Herfindahl-Hirschman index is used by the Department of Justice in the analysis of the competitive effects of mergers in specific markets [1]. It is calculated by squaring the market shares ($s_{i}$) of each firm in a market, summing the squares, and comparing this sum before and after the transaction under consideration.

\begin{equation}
H = 10^{4} \times \sum s_{i}^{2}
\end{equation}

By analogy, we propose a similar concentration measure for product portfolios, that may be applied to any slice of the portfolio along any path in the product hierarchy tree.

\begin{equation}
H = 10^{4} \times \sum \left(\frac{r_{i}}{\sum\limits ^n r_{i}}\right)^{2}
\end{equation}

$\left(\frac{r_{i}}{\sum\limits ^n r_{i}}\right)$ represents the share of revenue of each item within the slice of the product portfolio under consideration, the index $H$ is the sum of squared shares over the slice. By definition, $max(H) = 10^{4}$ and $min(H) = \frac{10^{4}}{n}$ (for a portfolio of $n$ items). See illustrative examples in tables 3 and 4.

\begin{table}[H]
\centering
\begin{tabular}{|l|l|l|l|} 
\hline
\textbf{Item}  & \textbf{Revenue ($r_{i}$)} & \textbf{Revenue Share} & \textbf{$10^{4} \times$Squared Revenue Share}  \\ 
\hline
Item 1         & \$180                      & 0.529                                                                     & 2803                                           \\ 
\hline
Item 2         & \$90                       & 0.265                                                                     & 701                                            \\ 
\hline
Item 3         & \$50                       & 0.147                                                                     & 216                                            \\ 
\hline
Item 4         & \$20                       & 0.059                                                                     & 35                                             \\ 
\hline
\textbf{Total} & \textbf{\$340}             & \textbf{1}                                                                & \textbf{H = 3754}                              \\
\hline
\end{tabular}
\caption{$H = 3754$ for a portfolio where half the revenue is concentrated in 1 item}
\end{table}

\begin{table}[H]
\centering
\begin{tabular}{|l|l|l|l|} 
\hline
\textbf{Item}  & \textbf{Revenue ($r_{i}$)} & \textbf{Revenue Share} & \textbf{$10^{4} \times$Squared Revenue Share}  \\ 
\hline
Item 1         & \$85                       & 0.25                                                                      & 625                                            \\ 
\hline
Item 2         & \$85                       & 0.25                                                                      & 625                                            \\ 
\hline
Item 3         & \$85                       & 0.25                                                                      & 625                                            \\ 
\hline
Item 4         & \$85                       & 0.25                                                                      & 625                                            \\ 
\hline
\textbf{Total} & \textbf{\$340}             & \textbf{1}                                                                & \textbf{H = 2500}                              \\
\hline
\end{tabular}
\caption{$H = 2500$ for a portfolio where revenue is evenly distributed across items}
\end{table}

Two graphical examples are also presented, one for a portfolio within a specific brand, where the largest item has a revenue share of less than 1\% and H = 18, the other for a portfolio within a specific brand, where the largest item has revenue share of greater than 15\% and H = 875. For the latter, with $t_{a} = 0.25$ only one item would be flagged as an A, which may not be commercially acceptable. To compensate for this concentration, we propose raising $t_{a}$.

\begin{figure}[H]

\includegraphics[scale=0.6]{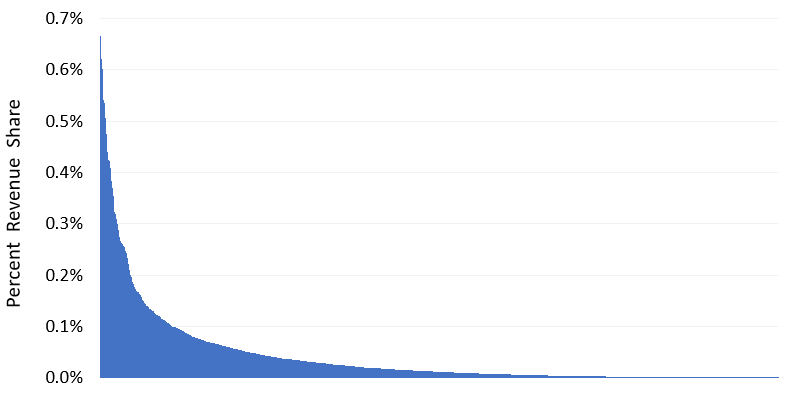}
\caption{H = 18 for a portfolio with the largest item having $<$1\% revenue share}

\end{figure}

\begin{figure}[H]
\includegraphics[scale=0.6]{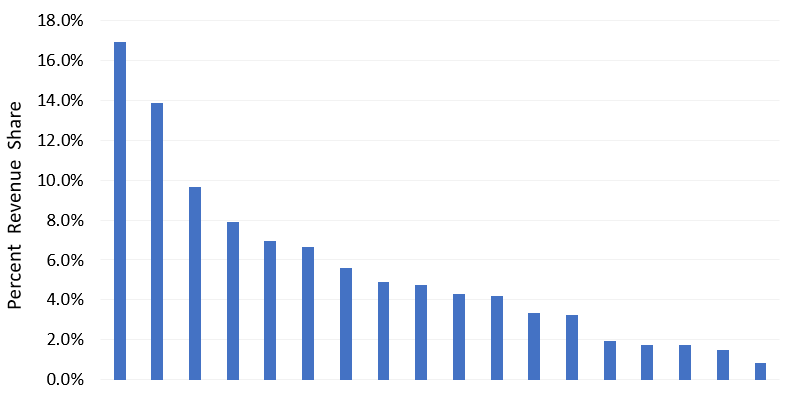}
\caption{H = 875 for a portfolio with the largest item having $>$15\% revenue share}
\end{figure}

\newpage

Some observations on $H$:

\begin{enumerate}
  \item For portfolios with the same item count ($n$), $H$ values may be directly compared

  \item For portfolios with different item counts and evenly distributed revenue $(\frac{1}{n})$, the ratio of $H$'s will equal the inverse ratio of item counts $(\frac{H_{2}}{H_{1}} = \frac{n_{1}}{n_{2}})$. For an ABCD algorithm, this special case is irrelevant as all items within each portfolio are identical and will receive the same class.
  
\end{enumerate}

We recommend reviewing the distribution of $H$ across the portfolio and simulating the impact on revenue and number of items classified as A's, to arrive at a commercially useful threshold and adjusted level of $t_{a}$. Using $H$ values to adjust thresholds for the ABCD algorithm is scalable as they are easily computed for large product portfolios. The thresholds for both $H$ and $t_{a}$ may be revisited annually for review.

\section*{Productivity} \label{s:sec4}
Productivity (average revenue per item across a group of items) is a widely used metric to characterize product portfolios, to track the evolution of the revenue distribution in portfolios over time and measure the progress of SKU (stock keeping unit) rationalization efforts. One can characterize the outcome of an ABCD analysis by measuring the productivity of each segment. For a portfolio of $n$ items, sorted by revenue high to low, productivity up to the $p^{th}$ item is

\begin{equation}
S_{p} = \frac{1}{p}\sum\limits ^{p} r_{i} \hspace{15 mm} \text{where } p \leq n
\end{equation}

with the following recurrence relation by definition

\begin{equation}
S_{p+1} = \frac{pS_{p} + r_{p+1} }{p+1}
\end{equation}

We now consider the case of the behavior of productivity when blending two portfolios. For a portfolio of $n$ items, sorted by revenue high to low, $S$ is by definition a non-increasing sequence. Now consider a second portfolio of $j$ items, with total revenue $J$ and productivity

\begin{equation}
S_{j} = \frac{J}{j}
\end{equation}

If we blend $all$ $j$ items from the second portfolio, with $p$ items from the first, the combined productivity ($T_{p}$) of the blended portfolio is  

\begin{equation}
T_{p} = \frac{pS_{p} + jJ }{p+j}
\end{equation}

where $j$ and $J$ are fixed for all values of $p$. Algebraically

\begin{equation}
T(p) = \frac{1}{p + j}\left (\int_{0}^{p}r(x)dx + jJ \right)
\end{equation}

with a maximum occuring when

\begin{equation}
\frac{dT(p)}{dp} = 0
\end{equation}

Solving Eq.[12] shows that maximum occurs when 

\begin{equation}
r(p^{*}) = T(p^{*})
\end{equation}

In other words, the combined productivity of the two portfolios is maximized at the $p^{th}$ item of the first, when the revenue of that item equals the combined productivity $T_{p}$ at that item. We illustrate this with an example in table 5. 

\begin{table}[H]
\centering
\begin{tabular}{|l|l|l|l|l|} 
\hline
\textbf{Item}   & \textbf{Revenue ($r_{i}$)} & \textbf{Item Count} & \textbf{Productivity~$(S_{p})$} & \textbf{Combined Productivity~$(T_{p})$}  \\ 
\hline
Item 1          & \$100                      & 1                                                             & \$100                         & \$40                                    \\ 
\hline
Item 2          & \$90                       & 2                                                             & \$95                          & \$50                                    \\ 
\hline
Item 3          & \$80                       & 3                                                             & \$90                          & \$55                                    \\ 
\hline
Item 4          & \$70                       & 4                                                             & \$85                          & \$57                                    \\ 
\hline
\textbf{Item 5} & \textbf{\$60}              & \textbf{5}                                                    & \textbf{\$80}                 & \textbf{\$58}                           \\ 
\hline
Item 6          & \$50                       & 6                                                             & \$75                          & \$57                                    \\ 
\hline
Item 7          & \$40                       & 7                                                             & \$70                          & \$55                                    \\ 
\hline
Item 8          & \$30                       & 8                                                             & \$65                          & \$53                                    \\ 
\hline
Item 9          & \$20                       & 9                                                             & \$60                          & \$50                                    \\
\hline
\end{tabular}
\end{table}

\begin{table}[H]
\begin{tabular}{|l|l|} 
\hline
Second portfolio combined revenue~$(J)$ & \$60  \\ 
\hline
Second portfolio item count~$(j)$       & 3     \\
\hline
\end{tabular}

\caption{The $1^{st}$ portfolio has 9 items and the $2^{nd}$ has 3 items with total revenue = \$60. When the entire second portfolio is combined with items from the $1^{st}$, sorted high to low, combined productivity is maximized at the $5^{th}$ item of the $1^{st}$ portfolio}

\end{table}

Now back to figure 2, if the value of $t_{a}$ in later passes is fixed, then it is possible that the total revenue \& item count, and therefore the productivity of A codes, from all the latter passes combined ($2^{nd}$ pass through the final pass) remain approximately constant. One can then solve for $t_{a}$ in the first pass to maximize the productivity of A items as a whole from the first stage of the algorithm, thus labeling additional items as A's.  

\begin{figure}[htb]
\includegraphics[scale=0.6]{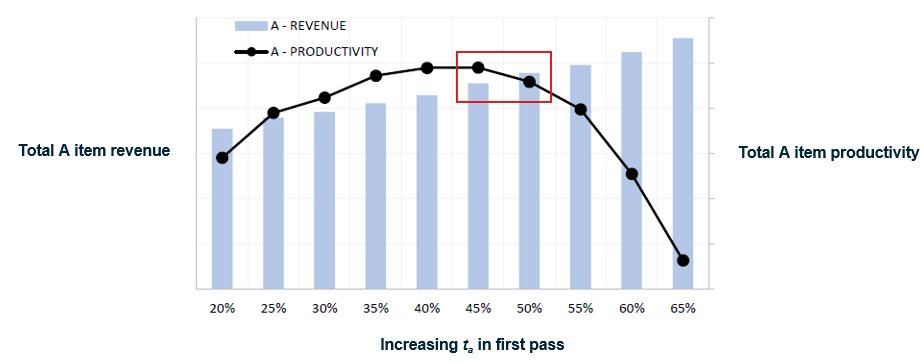}
\caption{A productivity maximum is achieved across A items, by varying $t_{a}$ in the first pass of the first $(A)$ stage of the algorithm (shown in figure 2). $t_{a}^{*} \approx 0.5$}
\end{figure}

The main point is that it may be possible to endogenously solve for $t_{a}$ in  special cases, and include additional items as A's, thus increasing commercial utility.

\section*{Summary} \label{s:sec5}
To conclude, it is our experience that a multi-stage, multi-pass ABCD algorithm, that provides for A items along different dimensions of the product hierarchy, leveraging a statistical concentration measure to calibrate thresholds, results in a more commercially viable classification for the supply chain. The key algorithmic features we have described can all be automated and scaled to large portfolios. Finally, the design and maintenance of a commercially viable and relevant ABCD segmentation algorithm is an iterative process, requiring constant refinement. Future directions of enquiry could include a joint optimization of the ABCD classification algorithm across other supply chain parameters, such as inventory policy and production prioritization.   

\section*{Acknowledgements} \label{s:sec6}
The authors would like to acknowledge the contributions of fellow Estee Lauder colleagues Raj Nakum, Chris Froah, Naresh Rajanna, Daniele Landi and Rogerio Pezutto. The authors also wish to thank Matt Olenich for diligently reviewing several early drafts.

\section*{Appendix} \label{s:sec7}

\begin{equation*}
T(p) = \frac{1}{p + j}\left (\int_{0}^{p}r(x)dx + jJ \right)
\end{equation*}

with a maximum occurring when

\begin{equation*}
\frac{dT(p)}{dp} = \frac{(p+j)r_{p} - 1 \cdot \left (\int_{0}^{p}r(x)dx + jJ \right)}{(p + j)^{2}} = 0
\end{equation*}

\begin{equation*}
r(p^{*}) = T(p^{*})
\end{equation*}

\end{document}